\documentclass[11pt]{article}
\usepackage{amsfonts}
\usepackage{color}
\usepackage{graphics}

\usepackage{cite}
\usepackage{latexsym}
\usepackage{amsmath}
\usepackage{amssymb}
\usepackage[dvips]{epsfig}
\usepackage{amscd}

 \usepackage{color}
\newtheorem{theorem}{Theorem}[section]
\newtheorem{remark}{Remark}[section]

\newtheorem{definition}{Definition}[section]
\newtheorem{lemma}[theorem]{Lemma}

\newtheorem{corollary}[theorem]{Corollary}

\newcommand{\n}{\rho}
\newcommand{\la}{\label}
\newcommand{\ti}{\tilde}

\newcommand{\lm}{\lambda}

\def\pf{{\it Proof.}  }
\renewcommand{\div}{ {\rm div }  }

\newcommand{\pa}{\partial}

\newcommand{\thatsall}{\hfill$\Box$}

\newcommand{\bt}{\begin{theorem}}
\newcommand{\bl}{\begin{lemma}}
\newcommand{\el}{\end{lemma}}
\newcommand{\et}{\end{theorem}}

\newcommand{\te}{\theta}

\newcommand{\de}{\delta}
\newcommand{\ve}{\varepsilon}

\newcommand{\bn}{\begin{eqnarray}}
\newcommand{\en}{\end{eqnarray}}
\newcommand{\bnn}{\begin{eqnarray*}}
\newcommand{\enn}{\end{eqnarray*}}

\newcommand{\bnnn}{\begin{eqnarray*}}
\newcommand{\ennn}{\end{eqnarray*}}
\newcommand{\ben}{\begin{enumerate}}
\newcommand{\een}{\end{enumerate}}

\newcommand{\ba}{\begin{aligned}}
\newcommand{\ea}{\end{aligned}}
\newcommand{\be}{\begin{equation}}
\newcommand{\ee}{\end{equation}}

\def\norm[#1]#2{\|#2\|_{#1}}

\def\lap{\triangle}

\def\lam{\lambda}

\def\rrr{\mathbb{R}^3}
\def\O{\Omega}

\def\r{\mathbb{R}}

\makeatletter
\@addtoreset{equation}{section}
\makeatother

\makeatletter
\@addtoreset{equation}{section}
\makeatother

\newcommand{\eps}{\varepsilon}
\newcommand{\Om}{\Omega}

\newcommand{\R}{\mathbb{R}}
\newcommand{\ds}{\displaystyle}

\newcommand{\al}{\alpha}
\newcommand{\ga}{\gamma}

\newcommand{\na}{\nabla}

\begin{document}
\title{Exponential Decay for Lions-Feireisl's Weak Solutions to the Barotropic Compressible Navier-Stokes Equations in   3D  Bounded Domains \thanks{ Partially supported by NNSFC 11501143, 11671027, and 11471321.}}
\date{}
\author{Yan-Fang Peng$^a$, Xiaoding Shi$^b$\thanks{
  Email addresses:  pyfang2005@sina.com (Y. Peng),   shixd@mail.buct.edu.cn (X. Shi)}  \\[3mm] $^{a}$School of Mathematical Sciences,   Guizhou Normal University,  \\Guiyang, 550001, P. R. China; \\
$^{b}$Department of Mathematics, School of Science, \\Beijing University of Chemical Technology, \\ Beijing, 100029, P. R. China }
\maketitle
\begin{abstract} For  barotropic compressible Navier-Stokes equations in three-dimensional (3D) bounded domains, we prove that any finite energy weak solution  obtained by  Lions [Mathematical topics in fluid mechanics, Vol. 2. Compressible models(1998)] and Feireisl-Novotn\'{y}-Petzeltov\'{a} [J. Math. Fluid Mech. 3(2001), 358-392]  decays exponentially to the equilibrium state. This result is established by both using the extra integrability of the density due to Lions and  constructing a suitable Lyapunov functional just under the framework of Lions-Feireisl's weak solutions.
\end{abstract}


Math Subject Classification: 35Q30; 76N10.

\section{Introduction}

  We consider the three-dimensional barotropic compressible Navier-Stokes equations    which read as follows:
\be\la{NS}
\begin{cases} \rho_t + \div(\rho u) = 0,\\
 (\rho u)_t + \div(\rho u\otimes u) + \nabla P = \mu\lap u + (\mu + \lam)\nabla \div u,
\end{cases}
\ee
where   $t\ge 0, x=(x_1,x_2,x_3)\in  \Omega\subset \rrr, \rho=\n(x,t)$ and $u=(u_1(x,t),u_2(x,t),u_3(x,t))$ represent, respectively,  the density and the velocity,  and the pressure $P$ is given by
\be\la{n2}
P(\rho) = A\rho^{\gamma}\, ( A>0,\, \ga>1).
\ee
 In the sequel, without loss of generality, we set $A= 1.$ The constant viscosity  coefficients  $\mu$ and  $\lambda$ satisfy  the  physical hypothesis:
\be\la{n3}
\mu>0,\quad 2\mu+3\lam\ge 0.
\ee

The initial conditions are imposed as
\begin{equation}\label{IV1}
\rho(x,t=0)=\rho_0,\quad\rho u(x,t=0)=m_0,
\end{equation}
together with the no-slip boundary conditions for the velocity
\begin{equation}\label{IV2}
u|_{\partial \Om}=0.
\end{equation}

A large number of literatures have been devoted to the large-time existence and behavior of solutions to compressible Navier-Stokes equations. On the one hand, for the existence of solutions, the one-dimensional problems have been studied extensively, see \cite{H,KS1} and the references therein. For the multidimensional case, the local existence and uniqueness of classical solutions were obtained in \cite{Na,Ser} in the absence of vacuum. The global classical solutions were first given by Matsumura-Nishida \cite{MN2} for the density strictly away from vacuum and the initial data close to a non-vacuum equilibrium. Later on, for discontinuous initial data, Hoff \cite{H3} showed  the existence of  global solutions  as limits of approximate solutions
corresponding to mollified initial data. Recently, Huang-Li-Xin \cite{HLX} first established global existence and uniqueness of  classical solutions with smooth initial data that are of small energy but possibly large oscillations and containing  vacuum states. As for the existence of weak solutions with large data, the major break-through is due to  Lions \cite{L2}, where for   three-dimensional case, the global weak solutions were obtained   under the condition  that  $ \ga\geq \frac95$ which was further relaxed  to  $\ga>\frac{3}{2}$ later  by  Feireisl-Novotny-Petzeltov\'{a} \cite{FNP}. Moreover, for the case that initial data has some symmetric properties, Jiang-Zhang \cite{JZ} proved that the equations possess   global weak solutions  for any $\ga>1.$

On the other hand, regarding   the large-time behavior of solutions for Navier-Stokes equations, Matsumura-Nishida \cite{MN1,MN2,MN3} proved the existence of global solutions  near a constant equilibrium state $(\rho_s,0)$ with $\n_s>0$ for the Cauchy problem in $\R^3,$ in the half space and exterior domains respectively. In particular, it was shown in \cite{MN2,MN3} that there exists a unique global classical solution $(\rho, u ) $ with $(\rho_0-\rho_s,m_0)$ sufficiently small in $H^3.$ Also, they showed that
$$
\|(\rho(\cdot,t)-\rho_s,  u(\cdot,t))\|_{L^\infty}=O(t^{-3/2}), \quad \|(\rho(\cdot,t)-\rho_s, u(\cdot,t))\|_{L^2}=O(t^{-3/4}),
$$
as $t\rightarrow \infty$ provided $(\rho_0-\rho_s,m_0)\in L^1.$  See also  \cite{HZ,PO,KS2,KK} and the references therein.
For the problem \eqref{NS}--\eqref{IV2}, Feiereisl-Petzeltov\'{a} \cite{FP1} showed that for any finite energy weak solution $(\rho,u),$ there exists  a stationary solution $(\rho_s,0)$ such that
$$
\ds \rho(\cdot,t)\rightarrow \rho_s~~\mbox{in} ~L^\ga(\Om),~\mbox{ess}\sup_{\tau>t}\int_\Om \rho(x,\tau)|u(x,\tau)|^2dx\rightarrow 0,~\mbox{as}~t\rightarrow \infty,
$$
where $\ga>\frac32,$ $\Om$ is not needed to be bounded and external force $\na F$  is independent of time $t.$ Moreover, Padula \cite{P} showed that in a bounded domain, the rest state is exponentially stable with respect to a large class of “weak” perturbations. More recently, Fang-Zi-Zhang \cite{FZZ} proved any finite energy weak solution to the problem \eqref{NS}--\eqref{IV2} without external force in bounded domains decays exponentially to the equilibrium state. However, the existence of  weak solutions considered by   \cite{P,FZZ} remains completely open for  large   data  since one of the basic assumptions in \cite{P,FZZ}  is that the density is bounded from above or below respectively uniformly in time. Indeed,    whether  Lions-Feireisl's finite energy weak solutions to the problem \eqref{NS}--\eqref{IV2} in a bounded domain whose existence is obtained by  \cite{FNP,L2}   decay  exponentially to  the equilibrium state or not remains open. In this paper, we will give a positive answer to this problem.

 Before  stating our main result, we    first introduce the definition of finite energy weak solutions.
\begin{definition}[\cite{FNP,L2}, Finite energy weak solutions]\label{Def1} A pair of functions $(\rho,u)$ will be termed as a finite energy weak solution of the problem \eqref{NS},\eqref{IV2} on $(0,\infty)\times \Om,$ if \\
$\bullet$~$\rho\geq 0,\,\,\rho\in L^{\infty}(0,\infty; L^\ga(\Om)),\, \, u\in L^2(0,\infty;H_0^1(\Om)).$\\
$\bullet$~The energy $$
E(t):=\int_{\Om}\left(\frac12 \rho |u|^2+\frac{1}{\ga-1}\rho^{\ga}\right)dx
$$  is locally integrable on $(0,\infty)$ and  for  any $0\le \psi(t)\in \mathcal{D}(0,\infty),$ it holds
\be\la{qj1} -\int_0^\infty \psi_tE(t)dt
+\int_0^\infty\psi \int_\Om \left[\mu|\na u|^2+(\lambda+\mu)(\div u)^2\right]dxdt\leq 0.
\ee \\
$\bullet$~Eqs. \eqref{NS} are satisfied in $\mathcal{D}'(0,\infty;\Om).$ Moreover, \eqref{NS}$_1$ holds in $\mathcal{D}'(0,\infty;\R^3)$ provided $(\rho,u)$ was prolonged to be zero on $\R^3\setminus {\Om}.$ \\
$\bullet$~\eqref{NS}$_1$ is satisfied in the sense of renormalized solutions, more precisely, the following equation
\begin{equation}\label{B}
b(\rho)_t+\div(b(\rho)u)+(b'(\rho)\rho-b(\rho))\div u=0
\end{equation}
holds in $\mathcal{D}'(0,\infty;\Om)$ for any $b\in C^1(\R)$ such that
\begin{equation*}\label{B1}
b'(z)\equiv 0,~ \mbox{for all} ~z\in \R~ \mbox{large enough},~ \mbox{say},~ |z|\geq M
\end{equation*}
where the constant $M$ may vary for different functions $b.$
\end{definition}
\begin{remark} A direct consequence of \eqref{qj1} is
\be\label{ES}
\ds \sup_{t\in (0,\infty)}E(t)+ \int_{0}^{\infty}\int_{\Om}\left[\mu|\nabla u|^2+(\lm+\mu)(\div u)^2\right]dxdt\leq E_0,
\ee with
\be\label{1.9}
 E_0:=\int_{\Om}\left(\frac{|m_0|^2}{2\rho_0}+\frac{1}{\ga-1}\rho_0^\ga\right) dx.
 \ee
\end{remark}

\begin{remark}\label{RB} It should be mentioned here that one can use the Lebesgue dominated convergence theorem to deduce that \eqref{B} will hold for any $b\in C^1(0,\infty)\cap C[0,\infty)$ satisfying
\begin{equation}\label{B2}
|b'(z)z|\leq c(z^\theta+z^{\frac{\ga}{2}}),\,\mbox{for all}\,~z>0\,\mbox{and a certain}\,~\theta\in (0,\frac{\ga}{2})
\end{equation}
provided $(\rho,u)$ is a finite energy weak solution in the sense of the above definition. In particular,  \eqref{B}  holds in $\mathcal{D}'(0,\infty;\Om)$  for $b(\n)=\n^\te$ with $\te\in (0,\ga/2].$
\end{remark}

Next, we state the following result concerning the existence of finite energy weak solution  to  problem \eqref{NS}--\eqref{IV2} due to Lions\cite{L2} and Feireisl-Novotny-Petzeltov\'{a} \cite{FNP}.

\begin{lemma}[\cite{FNP,L2}]\label{lm1.1} Assume $\Om\subset \R^3$ is a bounded domain of the class of $C^{2+\nu},\nu>0.$ Let $\ga>\frac{3}{2}$ and the initial data $(\rho_0, m_0)$ satisfy
$$ 0\le
\rho_0\in L^\ga(\Om), \,\,\frac{|m_0|^2}{\rho_0}\in L^1(\Om),
$$ with   $m_0=0$ almost everywhere on the set  $\{x\in\Om|\rho_0(x)=0\}.$
Then there exists  a finite energy weak solution $(\rho,u)$ of the problem \eqref{NS}--\eqref{IV2} satisfying for almost everywhere $t>0,$
\begin{equation}\label{1}
\ds \overline{ \rho}(t) =  \overline{ \rho_0} :=\n_s,
\end{equation} where (and in what follows),
$$\overline{f}:=\frac{1}{|\Om|}\int_\Om fdx $$
denotes the mean value of $f$ over $\Om,$
and for any $T>0,$ \be \la{hq1} \rho\in L^{ \ga+\te_0}(0,T;\Om), \ee   with some positive constant  $ \te_0\le -1+2\ga/3.$
\end{lemma}

\begin{remark}\label{RHO} As shown by Lions \cite{L2}, one can choose $\te_0=-1+2\ga/3.$
\end{remark}

Now we are in a position to state our main result as follows:
\begin{theorem}\label{Th} Assume that  the conditions of Lemma \ref{lm1.1} hold.   Then  there exist positive constants $C_1$ and $C_2$ both depending only on $\Om,\ga,\mu,\lm,\te_0,\n_s,$ and $ E_0$ such that $(\rho,u)$,  the  finite energy weak solution  to \eqref{NS}--\eqref{IV2}  whose existence is guaranteed by  Lemma \ref{lm1.1}, satisfies the following decay property:
\begin{equation}\label{1.10}
\ds\int_{\Om}(\rho|u|^2+G(\rho,\n_s))dx\leq C_2\exp\{-C_1t\}\,\,\mbox{a.e.}\,\,t>0,
\end{equation}
with
\be \la{gd1}G(\rho,\n_s):=\rho\int_{\rho_s}^\rho \frac{h^\ga-\rho_s^\ga}{h^2}dh.
\ee
\end{theorem}
\begin{remark}
It is worth  noticing that in Theorem \ref{Th}, our result holds for Lions-Feireisl's finite energy weak solutions whose existence is guaranteed by Lemma \ref{lm1.1}. Moreover, we do not attach any extra restriction on the weak solutions which indeed greatly improves those results of  \cite{P,FZZ} where the time-independent upper and/or lower bounds of density are essential in their analysis.\end{remark}

\begin{remark} After some small modifications, our method can be applied directly to other models, such as the
compressible magnetohydrodynamic flows in the barotropic case (see Appendix for the details), et al. \end{remark}

  We now make some comments on the analysis of this paper. To establish Theorem \ref{Th}, by combining  the energy inequality with the conservation of the mass (see \eqref{3.1}), the key issue is to discover new  decay estimates for $G(\n,\n_s)$ (see \eqref{gd1} for the definition).
    Compared with  \cite{P,FZZ} where  the time-independent upper and/or lower bounds of density play an essential role in their analysis, the main difficulties come  from  the fact that for  the finite energy weak solutions (see Definition \ref{Def1})  the    density is only  bounded time-independently on the $L^\infty(0,T;L^\ga(\Om))$-norm. To overcome these
    new difficulties,  we first observe that $G(\n,\n_s)$ can be bounded by $(\n^\ga-\n_s^\ga)(\n^\te-\n_s^\te)$ provided $\te>0$ (see \eqref{3.11}) and that for the finite energy weak solutions, the density has an additional integrality, that is, $\rho\in L^{\ga+\te_0}(0,T;\Om) $ for some positive  constant $\te_0$ (see \eqref{hq1}).
   Hence, to recover the decay estimate on the term      $(\n^\ga-\n_s^\ga)(\n^\te-\n_s^\te),$  we construct a  test function $\psi(t)\mathcal{B}([\rho^\theta]_\eps-\overline{[\rho^\theta]_\eps})$  for \eqref{NS}$_2$ with the aid of Bogovskii operator $\mathcal{B}$ (see Lemma \ref{lm2.6}) where we use the mollified functions $[\rho^\theta]_\eps$ and $\overline{[\rho^\theta]_\eps}$   due to  the lack of the integrality for $ \n^\te.$ Then after carefully using the commutator estimates (see \eqref{2.8}) and observing that (see \eqref{3.20})
$$
\int_0^\infty\psi  (\rho_s^\theta-\overline{\rho^\theta})
 \int_\Om (\rho^\ga-\rho_s^\ga) dxdt\geq 0,
 $$
 we obtain  the desired estimate on  $\int_0^\infty\psi \int_{\Om}G(\n,\n_s)dxdt$ (see \eqref{3.4b}) which is vital to get the decay estimate of the finite energy weak solutions.
   Finally,  observing that $(\n^\te-\n_s^\te)^2$ can be bounded by $G(\n,\n_s)$ provided $\te$ is suitably small (see \eqref{3.11}), we can build up a suitable  Lyapunov functional and then finish the proof of Theorem \ref{Th}.

This paper is organized as follows: in Sect.2, we establish some preliminary lemmas which will be needed in later analysis. In Sect.3, we are devoted to deriving some necessary estimates and finally prove Theorem \ref{Th}. Throughout the paper, $C$ denotes positive generic constant  depending only on $\Om,\ga,\mu,\lm,\te_0,\n_s,$ and $E_0$  which may vary in different cases. And we write $C(\al)$ to emphasize that $C$ depends on $\al.$

\section{ Preliminaries}

In this section, we recall some known facts and  elementary results which will be used later.

First, for $\eta$ as the standard mollifier in $\R^3$ and $f\in L^1_{loc}(\R^3),$ we set
$$
\eta_\eps(\cdot):=\frac{1}{\eps^3}\eta(\frac{\cdot}{\eps}),\quad [f]_\eps:=\eta_{\eps}\ast f.
$$

The following  properties of mollification are standard and can be found in \cite{A}.
\begin{lemma}[\cite{A}] \label{lm2.0} Let $f$ be a function which is defined on  $\R^3$ and vanishes identically outside a domain $\Om\subset\R^3.$

$(i)$~If $f\in L^1_{\rm loc}(\R^3),$ then $[f]_\eps\in C^\infty (\R^3).$

$(ii)$~If $f\in L^p(\Om)$ with $1\leq p<\infty,$ then $[f]_\eps\in L^p(\Om).$ Also
\be\label{2.1}
 \|[f]_\eps\|_{L^p(\Om)}\le \|f\|_{L^p(\Om)},\quad \lim_{\ve\rightarrow 0^+}\|[f]_\eps-f\|_{L^p(\Om)}=0.
\ee
\end{lemma}

Next, we state the commutator estimates which will play an important role in our further analysis.
\begin{lemma}[\cite{F,L1}]\label{lm2.4} Let $\Om\subset \R^N(N\geq2)$ be a domain and $f\in L^p(\Om), v\in [W^{1,q}(\Om)]^N$ be given functions with $1< p, q< \infty$ and $\frac1p+\frac1q\leq 1.$ Then for any compact $K\subset \Om,$
\bnn \begin{cases}
\|[\div(fv)]_\eps-\div([f]_\eps v)\|_{L^r(K)}\leq C(K) \|f\|_{L^p(\Om)}\|v\|_{W^{1,q}(\Om)},\\
\|[\div(fv)]_\eps-\div([f]_\eps v)\|_{L^r(K)}\rightarrow 0, ~\mbox{as}~\eps\rightarrow 0,\end{cases}
 \enn
provided $\eps>0$ is small enough and $\frac1r=\frac1p+\frac1q.$ In addition, if $\Om=\R^N,$ $K$ can be replaced by $\R^N.$
\end{lemma}

As a direct consequence of Lemma \ref{lm2.4}, we have

\begin{corollary}\label{co2.5}  Let $(\rho,u)$ be the solution of the problem \eqref{NS}--\eqref{IV2} as in Lemma \ref{lm1.1}. Then prolonging $(\rho,u)$ to zero in $\R^3\setminus {\Om},$ for $0<\theta\leq \frac{\ga}{2}$ and any $[\alpha,\beta]\subset (0,\infty),$ we have
\begin{equation}\label{2.9}
\partial_t([\rho^\theta]_\eps)+\div([\rho^\theta]_\eps u)=(1-\theta)[\rho^\theta \div u]_\eps+r_\eps, ~\mbox{a.e.~on}~ [\alpha,\beta]\times \Om,
\end{equation}
where $r_\eps:=\div([\rho^\theta]_\eps u)-[\div(\rho^\theta u)]_\eps.$
Moreover,
\be \label{2.8} r_\eps\rightarrow 0 \mbox{ in  }L^2(\alpha,\beta;L^q(\Om)),
\ee
for any $q\in [1,\frac{2\ga}{\ga+2\theta}].$
\end{corollary}
\pf First, taking $b(\rho)=\rho^\theta $ in \eqref{B}, we derive from \eqref{B2} that \begin{equation}\label{NS'}
(\rho^\theta)_t+\div(\rho^\theta u)+(\theta-1)\rho^\theta \div u=0, \mbox{ in}~\mathcal{D}'(0,\infty;\Om).
\end{equation}
Then, prolonging $(\rho,u)$ to zero in
$\R^3\setminus {\Om},$ we claim that
\begin{equation}\label{ES1}
\int_0^\infty \int_{\R^3}(\rho^\theta \varphi_t+\rho^\theta u\cdot\nabla \varphi)dxdt=(\theta-1)\int_0^\infty \int_{\R^3}\rho^\theta \varphi \div udxdt,
\end{equation}
for any $\varphi\in \mathcal{D}(0,\infty;\R^3).$ This in particular yields \eqref{2.9} provided $\eps>0$ is small enough. Since $u\in L^2(0,\infty;H_0^1(\Om)),$ Lemma \ref{lm2.4} shows that $r_\eps$ is bounded in $L^2(\alpha,\beta;L^{\frac{2\ga}{\ga+2\theta}}(\Om))$ uniformly in $\ve.$ By Lemma \ref{lm2.4} and  the Lebesgue dominated convergence theorem,  we have  $r_\eps\rightarrow 0$ in $L^2(\alpha,\beta;L^{\frac{2\ga}{\ga+2\theta}}(\Om))$ which together with the boundedness of the domain $\Om$ gives \eqref{2.8}.

Finally, it only remains to prove \eqref{ES1}. To this end, we mainly extract some ideas from \cite[Lemma 3.3]{FNP} and take a sequence of functions $\phi_m\in \mathcal{D}(\Om)$ satisfying
\begin{equation}\label{PHI}
\left\{%
\begin{array}{ll}
0\leq \phi_m\leq 1; \phi_m=1, ~\mbox{for all}~x ~\mbox{such that dist}(x,\pa \Om)\geq \frac1m,\\
|\nabla \phi_m(x)|\leq 2m,~ \mbox{for all}~ x\in \Om.
\end{array}
\right.
\end{equation}
Then
\be \la{KJ1}
 \ds \int_0^\infty\int_{\R^3}\rho^\theta \varphi_tdxdt=\int_0^\infty\int_{\Om}\rho^\theta(\phi_m\varphi)_t dx dt+\int_0^\infty\int_{\Om}\rho^\theta(1-\phi_m)\varphi_tdxdt
 \ee
 and
\be\label{3.4}\ba&\int_0^\infty\int_{\R^3}\rho^\theta u \cdot\na \varphi dxdt
 \\&= \int_0^\infty\int_{\Om}\rho^\theta u\cdot \na(\phi_m\varphi)dxdt+\int_0^\infty\int_{\Om}(1-\phi_m)\rho^\theta u\cdot\na \varphi dxdt\\
 & \quad-\int_0^\infty\int_{\Om}\varphi\rho^\theta u\cdot \na \phi_mdxdt.
 \ea\ee

Moreover, by \eqref{NS'}, we have
\bnn\ba&
\int_0^\infty\int_{\Om}\rho^\theta(\phi_m\varphi)_tdxdt+\int_0^\infty\int_{\Om} \rho^\theta u \cdot\na(\phi_m\varphi)dxdt\\&=
(\theta-1)\int_0^\infty\int_{\Om}\rho^\theta (\phi_m\varphi)\div u dxdt,
\ea\enn
which together with \eqref{ES1}, \eqref{KJ1}, and  \eqref{3.4}   gives
  \be\la{pj1}\ba
 & \int_0^\infty\int_{\R^3}(\rho^\theta \varphi_t+\rho^\theta u \cdot\na \varphi)dxdt\\
 & =(\theta-1)\int_0^\infty\int_{\Om}\rho^\theta(\phi_m\varphi) \div u dxdt+ \int_0^\infty\int_{\Om}\rho^\theta (1-\phi_m)\varphi_tdxdt\\
 & \quad+\int_0^\infty\int_{\Om}(1-\phi_m)\rho^\theta u\cdot\na \varphi dxdt-\int_0^\infty\int_{\Om}\varphi\rho^\theta u\cdot\na \phi_mdxdt.
\ea\ee
Then, on the one hand,  it follows from \eqref{ES}, \eqref{PHI} and the Lebesgue dominated convergence theorem that as $m\rightarrow \infty,$
\be \la{pj2}
\int_0^\infty\int_{\Om}\rho^\theta(\phi_m\varphi) \div u dxdt\rightarrow \int_0^\infty\int_{\Om}\rho^\theta\varphi \div u dxdt,
\ee
\be\label{2.11}
\int_0^\infty\int_{\Om}\rho^\theta(1-\phi_m)\varphi_t dxdt\rightarrow 0,
\ee
and
\be\label{2.12}
\int_0^\infty\int_{\Om}\rho^\theta u(1-\phi_m)\varphi_t dxdt\rightarrow 0,
\ee
due to $0< \theta\leq \frac{\ga}{2}.$

On the other hand, it follows from \eqref{PHI} and the Lebesgue dominated convergence theorem that
\be \label{2.13}\ba
&\left|\int_0^\infty\int_{\Om}\varphi \rho^\theta u\cdot\na \phi_mdxdt\right|\\ &\le
C\left(\int_0^\infty\int_{\Om}1_{(dist(x,\pa \Om)\le 1/m)} \varphi^2 \rho^{2\theta} dxdt\right)^{1/2}\left( \int_0^\infty\int_{\Om} \frac{|u|^2}{{ dist(x,\pa \Om)}^{2} }dxdt\right)^{1/2}\\ &\le
C\left(\int_0^\infty\int_{\Om}1_{(dist(x,\pa \Om)\le 1/m)} \varphi^2 \rho^{2\theta} dxdt\right)^{1/2}\left( \int_0^\infty\int_{\Om}  |\na u|^2  dxdt\right)^{1/2}\\&\rightarrow 0,~\mbox{as}~ m\rightarrow \infty,\ea\ee
where in the second inequality we have used  Hardy's inequality and $u\in L^2(0,\infty;H_0^1(\Om)).$ Taking $m\rightarrow \infty$ in \eqref{pj1} and using \eqref{pj2}--\eqref{2.13} leads to \eqref{ES1}, then we finish the proof of Corollary  \ref{co2.5}. \thatsall

Finally, let $\Om$ be a bounded Lipschitz domain in $\R^3.$ We consider an auxiliary problem
 \begin{equation}\label{OP}
\div~ v=f,~~v|_{\pa \Om}=0.
\end{equation}
\begin{lemma} [\cite{B,GA}] \label{lm2.6} For problem \eqref{OP},  there exists a linear operator $\mathcal{B}=[\mathcal{B}_1,\mathcal{B}_2,\mathcal{B}_3]$ enjoying the properties:

$\bullet$~$\mathcal{B}$ is a bounded linear operator from $\{f\in L^p(\Om)|\int_\Om fdx=0\}$ into $[W_0^{1,p}(\Om)]^3,$ that is,
$$
\|\mathcal{B}[f]\|_{W_0^{1,p}(\Om)}\leq C(p,\Om)\|f\|_{L^p(\Om)},\,\mbox{for any }p\in(1,\infty).
$$

$\bullet$~The function $v=\mathcal{B}[f]$ solves  the problem \eqref{OP}.

$\bullet$~If, moreover, $f $ can be written in the form $f=\div g$ with $g\in [L^r(\Om)]^3$ and $g\cdot \vec{n}|_{\pa \Om}=0,$  then
$$
\|\mathcal{B}[f]\|_{L^r(\Om)}\leq C( r,\Om)\|g\|_{L^r(\Om)}, ~\mbox{for any}\,\,r\in (1,\infty).
$$
\end{lemma}
\begin{remark}
The  operator $\mathcal{B}$ was first constructed by the Bogovskii \cite{B}. A complete proof of the above mentioned properties may be found in Galdi \cite[Theorem 3.3]{GA}.
\end{remark}

\section{Proof of Theorem \ref{Th}}

First, recalling that $\rho_s=\overline{\rho_0 }$ is a positive constant, we deduce from \eqref{1} that for  any $0\le \psi(t)\in \mathcal{D}(0,\infty),$
$$
\ds \int_0^\infty \psi_t\int_\Om\left(-\frac{\ga}{\ga-1}\rho \rho_s^{\ga-1}+\rho_s^\ga\right)dxdt=0,
$$
which together with \eqref{qj1}  gives
 \begin{equation}\label{3.1}\ba&
 -\int_0^\infty \psi_t\int_\Om \left(\frac12 \rho|u|^2+G(\rho,\n_s)\right)dxdt\\&+\int_0^\infty \psi\int_\Om \left[\mu |\na u|^2+(\lm+\mu)(\div u)^2\right]dxdt\leq 0.\ea
 \end{equation}

 Next, noticing that for $ \theta:= \min\{1,\te_0,{\ga}\}/4,$ there exist  positive constants $C_0$ and $\ti C_0$ both depending only on $\ga,\te_0$ and $\n_s$ such that for any $\n\ge 0,$  \be\label{3.11} \ti C_0( \rho^\theta -  \rho_s^\theta )^2\le C_0  G(\n,\n_s)  \le    (\n^\ga-\n_s^\ga)( \rho^\theta - \rho_s^\theta). \ee
Then we claim that there exists some  constant $C>0$ depending only on $\Om,\ga,\mu,\lm,\te_0,$ and $ E_0$ such that for any $0\le \psi(t)\in \mathcal{D}(0,\infty),$
\be\label{3.4b}
\ba&  \frac{C_0}{2}\int_0^\infty\psi\int_{\Om}G(\n,\n_s) dxdt+\int_0^\infty \psi_t\int_{\Om}\rho u\cdot\mathcal{B}( \rho^\theta -\overline{ \rho^\theta })dxdt\\ &\le C\int_0^\infty \psi \|\na u\|_{L^2}^2dt.
\ea\end{equation}

Adding \eqref{3.4b} multiplied by a suitably small constant $\delta>0$ which will be determined later to \eqref{3.1} gives
 \begin{equation}\label{3.18s}
 \ba&
 -\int_0^\infty \psi_t V_\de(t)   dt +\int_0^\infty \psi  W_\de (t) dt \leq 0,\ea
 \end{equation}
where
 \bnn\label{3.24} V_\de(t):=\int_\Om\left(\frac12 \rho|u|^2+G(\rho,\n_s)-\de \rho u\cdot\mathcal{B}( \rho^\theta -\overline{ \rho^\theta }) \right)dx,\enn and \bnn W_\de(t):=\int_\Om\left((\mu-C\de) |\na u|^2+(\lm+\mu)(\div u)^2+\frac{C_0\de}{2}G(\n,\n_s) \right)dx.
 \enn

Moreover, it follows from   Lemma \ref{lm2.6}  and \eqref{ES} that
\bnn
\ba
  & \left|\int_{\Om} \rho u\cdot\mathcal{B}(\rho^\theta-\overline{\rho^\theta})dx\right|\\
  &
  \leq \ds  \|\sqrt{\rho}u\|_{L^2(\Om)}\|\sqrt{\rho}\|_{L^3(\Om)}\|\mathcal{B}(\rho^\theta-\overline{\rho^\theta})\|_{L^6(\Om)}\\
  &\leq \ds C \|\sqrt{\rho}u\|_{L^2(\Om)}\|\na \mathcal{B}(\rho^\theta-\overline{\rho^\theta})\|_{L^2(\Om)}\\
  &\leq \ds  \int_\Om \rho|u|^2dx+C  \|\rho^\theta-\overline{\rho^\theta}\|^2_{L^2(\Om)}\\
  &\leq \ds \int_\Om \rho|u|^2dx+C  \|\rho^\theta-\rho_s^\theta\|^2_{L^2(\Om)}+C\|\overline{\rho^\theta}-\overline{\rho^\theta_s}\|^2_{L^2(\Om)} \\
  &\leq      \int_\Om \rho|u|^2dx+C  \int_\Om G(\rho,\n_s)dx
  \ea\enn
 where in the last inequality we have used \eqref{3.11}. Then after choosing  suitably small ${\de_0}>0$ which depends only on $\Om,\ga,\mu,\lm,\te_0, \n_s,$ and $ E_0,$  we get
\begin{equation}\label{3.17}
\frac14 \int_\Om \left( \rho |u|^2+G(\rho,\n_s)\right)dx\le V_{\de_0}(t)\leq 2 \int_\Om \left( \rho |u|^2+G(\rho,\n_s)\right)dx\ee
and
\be\label{3.27} W_{\de_0}(t)\ge \frac14\int_\Om\left( \mu  |\na u|^2 + {C_0{\de_0}} G(\n,\n_s) \right)dx.
\end{equation}

Note that
\bnn\ba \int_\Om \n |u|^2dx \le \|\n\|_{L^{\frac32}(\Om)}\|u\|_{L^6(\Om)}^2 \le C\|\na u\|_{L^2(\Om)}^2,\ea
\enn
which together with \eqref{3.17} and \eqref{3.27} implies for almost everywhere $t\in(0,\infty),$
\bnn  V_{\de_0}(t)\le C_1 W_{\de_0}(t), \enn with some constant $C_1>0$ depending   on $\Om,\ga,\mu,\lm,\te_0, \n_s,$ and $ E_0.$
Putting this into \eqref{3.18s} yields that for any $0\leq\psi\in \mathcal{D}(0,\infty),$
\begin{equation}\label{3.9}
-\int_0^\infty \psi_t V_{\de_0}(t) dt
  +C_1\int_0^\infty \psi V_{\de_0}(t) dt\leq 0.
\end{equation}

Let $[a,b]$ be any compact subset of $(0,\infty).$  Taking $\psi(t)=\eta_\eps(t-\cdot)$ in \eqref{3.9} gives
$$
\partial_t[V_{\de_0}]_\eps+C_1[V_{\de_0}]_\eps\leq 0,~\mbox{a.e.}~ t\in [a,b],
$$
provided $\eps>0$ is small enough.
 This implies
$$
[V_{\de_0}]_\eps(t)\leq [V_{\de_0}]_\eps(s)\exp\{-C_1(t-s)\}
$$
for a.e. $0<s<t<\infty.$   Since $V_{\de_0}(t)\in L_{loc}^\infty(0,\infty),$ letting $\eps\rightarrow 0,$ we have
\bnn
V_{\de_0}(t)\leq V_{\de_0}(s)\exp\{-C_1(t-s)\},
\enn
which together with \eqref{3.17}, \eqref{ES} and \eqref{1} yields  \eqref{1.10}.

In the end,  it only remains to prove \eqref{3.4b}.
Indeed, first,  set $$
 \Phi(x,t)=\psi(t)\mathcal{B}([\rho^\theta]_\eps-\overline{[\rho^\theta]_\eps}).
 $$
 Since $\rho\in L^{\ga+\theta_0}(0,T;\Om),$
by Corollary \ref{co2.5}, we can use $\Phi(x,t)$ as a test function for \eqref{NS}$_2$ to get
 \be\la{s3.4} \ba  &  \int_0^\infty\psi \int_{\Om}(P(\rho)-P(\rho_s))([\rho^\theta]_\eps -\overline{[\rho^\theta]_\eps})dxdt +\int_0^\infty \psi_t\int_{\Om}\rho u\cdot\mathcal{B}( \rho^\theta -\overline{ \rho^\theta })dxdt \\
 &=-\int_0^\infty \psi_t\int_{\Om}\rho u\cdot\mathcal{B}\left([\rho^\theta]_\eps-\overline{[\rho^\theta]_\eps}-( \rho^\theta -\overline{ \rho^\theta })\right)dxdt \\&\quad-\int_0^\infty \psi \int_{\Om}\rho u \cdot \mathcal{B}\left(([\rho^\theta]_\eps)_t-\overline{([\rho^\theta]_\eps)_t}\right)dxdt \\
 &\quad-\int_0^\infty\psi\int_{\Om}\rho u\otimes u:\na \mathcal{B} \left([\rho^\theta]_\eps-\overline{[\rho^\theta]_\eps}\right)dxdt\\
 & \quad+\int_0^\infty \psi\int_\Om \left(\mu \na u:\na \mathcal{B}([\rho^\theta]_\eps-\overline{[\rho^\theta]_\eps})
  +(\lm+\mu)  \div u([\rho^\theta]_\eps-\overline{[\rho^\theta]_\eps})\right)dxdt\\
 &:\ds =-\int_0^\infty \psi_t I_1dt+\int_0^\infty \psi\sum_{i=2}^4I_i dt. \ea\ee

Using Lemmas \ref{lm2.0}, \ref{lm2.6} and \eqref{ES}, we estimate each  $I_i(i=1,\cdots,4)$   as follows:
\be \la{s3.5}\ba
 |I_1|
& \leq \ds \|\sqrt{\rho}u\|_{L^2(\Om)}\|\sqrt{\rho}\|_{L^3(\Om)} \|\mathcal{B}([\rho^\theta]_\eps-\overline{[\rho^\theta]_\eps}-( \rho^\theta -\overline{ \rho^\theta }))\|_{L^6(\Om)}\\
& \leq \ds C\|\sqrt{\rho}u\|_{L^2(\Om)} \|\na \mathcal{B}([\rho^\theta]_\eps-\overline{[\rho^\theta]_\eps}-( \rho^\theta -\overline{ \rho^\theta }))\|_{L^2(\Om)}  \\
& \leq \ds C  \|[\rho^\theta]_\eps-\rho^\theta\|_{L^2(\Om)},\ea\ee and
\be\label{s3.13}\ba
|I_3|
& \leq \ds \|\rho\|_{L^\ga(\Om)}\|u\|^2_{L^6(\Om)}\|\na \mathcal{B}([\rho^\theta]_\eps-\overline{[\rho^\theta]_\eps})\|_{L^{\frac{3\ga}{2\ga-3}}(\Om)}\\
& \leq \ds C\|\na u\|^2_{L^2(\Om)} \|[\rho^\theta]_\eps-\overline{[\rho^\theta]_\eps}\|_{L^{\frac{3\ga}{2\ga-3}}(\Om)}\\
& \leq C\|\na u\|^2_{L^2(\Om)} \|\rho^\theta\|_{L^{\frac{3\ga}{2\ga-3}}(\Om)}\\
& \leq C\|\na u\|^2_{L^2(\Om)},
\ea\ee
due to  $\theta\leq \frac{2\ga-3}{3}.$

Furthermore, direct computation shows that for $1\leq p\leq 1+\frac{\ga}{\theta},$
\bnn\ba
   &\ds \|[\rho^\theta]_\eps-\overline{[\rho^\theta]_\eps}  \|_{L^p(\Om)} \\
  &\leq \ds \|[\rho^\theta]_\eps-\rho^\theta\|_{L^p(\Om)}+\|\rho^\theta-\rho_s^\theta\|_{L^p(\Om)}+\|\rho_s^\theta-\overline{\rho^\theta}\|_{L^p(\Om)}
 +\|\overline{\rho^\theta}-\overline{[\rho^\theta]_\eps}\|_{L^p(\Om)}\\
  &\leq \ds C\|\rho^\theta-\rho_s^\theta\|_{L^p(\Om)}+C\|[\rho^\theta]_\eps -\rho^\theta\|_{L^p(\Om)},
\ea\enn
which gives
\be \ba
 |I_4|
& \leq C \|\na u\|_{L^2(\Om)} \|  [\rho^\theta]_\eps-\overline{[\rho^\theta]_\eps} \|_{L^2(\Om)}\\
& \leq \ds C\|\na u\|_{L^2(\Om)}\| \rho^\theta - \rho_s^\theta \|_{L^2(\Om)} +C\|\na u\|_{L^2(\Om)}\|[\rho^\theta]_\eps-  \rho^\theta \|_{L^2(\Om)}\\
& \leq  \frac{1}{4} \int_\O (\n^\ga-\n_s^\ga)( \rho^\theta - \rho_s^\theta)dx  +C\|\na u\|_{L^2(\Om)}^2 + C\|[\rho^\theta]_\eps-  \rho^\theta \|^2_{L^2(\Om)},
\ea\ee
where in the last inequality we have used \eqref{3.11}.

 As for $I_2,$ by \eqref{2.9}, we have
 \be \label{3.12}\ba
  I_2  & =  -(1-\theta) \int_{\Om}\rho u\cdot\mathcal{B}([\rho^\theta \div u]_\eps-\overline{[\rho^\theta \div u]_\eps})dx  \\
  &\quad +  \int_{\Om}\rho u \cdot\mathcal{B}(\div([\rho^\theta]_\eps u))dx-  \int_{\Om}\rho u\cdot\mathcal{B}(r_\eps-\overline{r_\eps})dx   :  = \sum_{i=1}^3I_2^i.
\ea
\ee
It follows from Lemmas \ref{lm2.0} and \ref{lm2.6} that for $\tilde{\ga}:=\min\{\ga,5\},$
\be\ba\label{3.13}
 |I^1_2|
& \leq C\|\rho\|_{L^{\tilde{\ga}}(\Om)}\|u\|_{L^6(\Om)}\|\mathcal{B} ([\rho^\theta \div u]_\eps-\overline{[\rho^\theta \div u]_\eps})\|_{L^{\frac{6\tilde{\ga}}{5\tilde{\ga}-6}}(\Om)}\\
& \leq \ds C\|\na u\|_{L^2(\Om)}\|\na  \mathcal{B}([\rho^\theta \div u]_\eps-\overline{[\rho^\theta \div u]_\eps})\|_{L^{\frac{6\tilde{\ga}}{7\tilde{\ga}-6}}(\Om)}\\
& \leq \ds C\|\na u\|_{L^2(\Om)} \|[\rho^\theta \div u]_\eps-\overline{[\rho^\theta \div u]_\eps}\|_{L^{\frac{6\tilde{\ga}}{7\tilde{\ga}-6}}(\Om)}\\
& \leq \ds C\|\na u\|_{L^2(\Om)} \|\rho^\theta \div u\|_{L^{\frac{6\tilde{\ga}}{7\tilde{\ga}-6}}(\Om)}\\
& \leq \ds C\|\na u\|_{L^2(\Om)}\|\rho^\theta\|_{L^{\frac{3\tilde{\ga}} {2\tilde{\ga}-3}}(\Om)}\|\na u\|_{L^2(\Om)}\\
& \leq \ds C\|\na u\|^2_{L^2(\Om)},
\ea\ee
where in the last inequality we have used \eqref{ES} and  $\theta\leq \frac{\ga(2\tilde{\ga}-3)}{3\tilde{\ga}}.$

Similarly,
\be\la{s3.14}\ba
  |I_2^2| &\leq C\|\rho\|_{L^\ga(\Om)} \|u\|_{L^6(\Om)} \|\mathcal{B}(\div([\rho^\theta]_\eps u)\|_{L^{\frac{6\ga}{5\ga-6}}(\Om)}\\
& \leq C\|\na u\|_{L^2(\Om)} \|[\rho^\theta]_\eps u\|_{L^{\frac{6\ga}{5\ga-6}}(\Om)}\\
& \leq C\|\na u\|_{L^2(\Om)} \|u\|_{L^6(\Om)}\|\rho^\theta\|_{L^{\frac{3\ga}{2\ga-3}}(\Om)} \\
& \leq C\|\na u\|^2_{L^2(\Om)},
\ea\ee
due to $\theta\leq \frac{2\ga-3}{3}.$ Moreover,
\bnn\ba
 |I_2^3| &\leq  \|\rho\|_{L^{\tilde{\ga}}(\Om)} \|u\|_{L^6(\Om)} \|\mathcal{B}(r_\eps-\overline{r_\eps})\|_{L^\frac{5\tilde{\ga}}{5\tilde{\ga}-6}(\Om)}\\
& \leq C\|\na u\|_{L^2(\Om)}\|\na \mathcal{B}(r_\eps-\overline{r_\eps})\|_{L^{\frac{6\tilde{\ga}}{7\tilde{\ga}-6}}(\Om)}\\
& \leq C\|\na u\|_{L^2(\Om)}\|r_\eps-\overline{r_\eps}\|_{L^{\frac{6\tilde{\ga}}{7\tilde{\ga}-6}}(\Om)}\\
& \leq C \|\na u\|_{L^2(\Om)}\|r_\eps\|_{L^{\frac{6\tilde{\ga}}{7\tilde{\ga}-6}}(\Om)},
\ea\enn
which together with \eqref{ES} and \eqref{2.8}  leads to
\begin{equation}\label{s3.15}\ba
\int_0^\infty \psi |I_2^3|dt &\le C\left(\int_0^\infty \psi \|\na u\|_{L^2(\Om)}^2dt\right)^{1/2}\left(\int_0^\infty \psi\|r_\eps\|^2_{L^{\frac{6\tilde{\ga}}{7\tilde{\ga}-6}}(\Om)}dt\right)^{1/2} \\&\rightarrow 0,~\mbox{as}~ \eps\rightarrow 0,\ea
\end{equation}
due to $\theta\leq \frac{\ga(2\tilde{\ga}-3)}{3\tilde{\ga}}.$

 Next, for the first term on the left-hand side of \eqref{s3.4},  we have
\be\la{s3.12}\ba
& \int_0^\infty\psi\int_{\Om}(P(\rho)-P(\rho_s))([\rho^\theta]_\eps-\overline{[\rho^\theta]_\eps})dxdt\\
& =\int_0^\infty\psi\int_{\Om}(\rho^\ga-\rho_s^\ga)\left[([\rho^\theta]_\eps-\rho^\theta)+(\rho^\theta-\rho_s^\theta)+(\rho_s^\theta-\overline{\rho^\theta})
+(\overline{\rho^\theta}-\overline{[\rho^\theta]_\eps})\right]dxdt\\
& \geq \int_0^\infty\psi\int_\Om (\rho^\ga-\rho_s^\ga)(\rho^\theta-\rho_s^\theta) dxdt+\int_0^\infty\psi\int_\Om (\rho^\ga-\rho_s^\ga)([\rho^\theta]_\eps-\rho^\theta)dxdt \\
& \quad+\int_0^\infty\psi (\overline{\rho^\theta}-\overline{[\rho^\theta]_\eps})\int_\Om (\rho^\ga-\rho_s^\ga)dxdt\\
& : =\int_0^\infty\psi\int_\Om (\rho^\ga-\rho_s^\ga)(\rho^\theta-\rho_s^\theta) dxdt +J,
\ea\ee
where in the second inequality we have  used  the following inequality:
\be \label{3.20}
\int_0^\infty\psi  (\rho_s^\theta-\overline{\rho^\theta})
 \int_\Om (\rho^\ga-\rho_s^\ga) dxdt\geq 0,
\ee
due to the following simple fact:
\bnn
\overline{\rho^\ga}\ge \overline{\rho}^\ga=\rho_s^\ga  ,\quad\rho_s^\theta=\overline{\rho}^\te\ge \overline{\rho^\theta}.
\enn

As for $J,$ since $\theta\leq\min\{\ga/2,\theta_0\},$ it follows from H\"older's inequality and \eqref{hq1} that
\bnn\ba |J|&\leq    \int_0^\infty \psi \|\rho^\ga-\rho_s^\ga\|_{L^{\frac{\ga+\theta}{\ga}}(\Om)} \|[\rho^\theta]_\eps-\rho^\theta\|_{L^{\frac{\ga+\theta}{\theta}}(\Om)}dt\\&\quad+C\int_0^\infty \psi \| [\rho^\theta]_\eps-\rho^\theta\|_{L^1(\Om)}dt\\
& \leq C(\psi)\left(\int_0^\infty \psi\|[\rho^\theta]_\eps-\rho^\theta\| ^{\frac{\ga+\theta}{\theta}}_{L^{\frac{\ga+\theta}{\theta}}(\Om)}dt\right) ^{\frac{\theta}{\ga+\theta}}+ C\int_0^\infty \psi \| [\rho^\theta]_\eps-\rho^\theta\|_{L^2(\Om)}dt,
\ea\enn
which together with  \eqref{s3.4}-\eqref{s3.14}  and  \eqref{s3.12} leads to
\be\label{3.4a}
\ba
&  \frac12\int_0^\infty\psi\int_{\Om}(\rho^\ga-\rho_s^\ga)(\rho^\theta-\rho_s^\theta) dxdt+\int_0^\infty \psi_t\int_{\Om}\rho u\cdot\mathcal{B}( \rho^\theta -\overline{ \rho^\theta })dxdt\\ &\le C\int_0^\infty \psi \|\na u\|_{L^2}^2dt+ C\int_0^\infty(|\psi_t| +\psi)\|[\rho^\theta]_\eps-  \rho^\theta \|_{L^2(\Om)}dt \\&\quad+C(\psi)\left(\int_0^\infty \psi\|[\rho^\theta]_\eps-\rho^\theta\| ^{\frac{\ga+\theta}{\theta}}_{L^{\frac{\ga+\theta}{\theta}}(\Om)}dt\right) ^{\frac{\theta}{\ga+\theta}}+C\int_0^\infty\psi |I_2^3|dt\\&\quad+C\int_0^\infty\psi \|[\rho^\theta]_\eps-  \rho^\theta \|^2_{L^2(\Om)}dt\\&:=C\int_0^\infty \psi \|\na u\|_{L^2}^2dt+M(\ve).
\ea\end{equation}

Finally, it follows from \eqref{3.4}, \eqref{hq1}, Lemma \ref{lm2.0}$(ii)$ and the Lebesgue dominated convergence theorem that \bnn \lim_{\ve \rightarrow 0 }M(\ve)= 0, \enn which together with \eqref{3.4a} and \eqref{3.11}  gives \eqref{3.4b}. The proof of  Theorem \ref{Th} is completed. \thatsall

  \appendix
  \renewcommand{\appendixname}{Appendix~\Alph{section}}

  \section{Appendix}
 In this section,  we will show how to apply our method to study other models. As an example, we consider the equations of three-dimensional
compressible magnetohydrodynamic flows in the barotropic case as follows(\cite{CH,KL,LL}):
\be\la{MHD}
\begin{cases} \rho_t + \div(\rho u) = 0,\\
 (\rho u)_t + \div(\rho u\otimes u) + \nabla P = (\nabla \times H)\times H+\mu\lap u + (\mu + \lam)\nabla \div u,\\
 H_t-\nabla \times(u\times H)=-\nabla \times (\nu \nabla \times H),\quad \div H=0,
\end{cases}
\ee
where $\n$ denotes the density, $u\in\r^3$ the velocity, $H\in\r^3$ the magnetic field,
$P(\n)=A\n^\ga$ the pressure with constant $A>0$ and the adiabatic exponent $\ga>1$;
the viscosity coefficients of the flow satisfy $2\mu+3\lambda> 0$ and $\mu>0$; $\nu>0$ is
the magnetic diffusivity acting as a magnetic diffusion coefficient of the magnetic
field, and all these kinetic coefficients and the magnetic diffusivity are independent
of the magnitude and direction of the magnetic field. We impose  the following initial-boundary conditions on \eqref{MHD}:
 \be\la{IC}
\begin{cases}
\rho(x,0)=\rho_0(x)\in L^\ga(\Om),\quad \rho_0(x)\geq 0,\\
(\rho u)(x,0)=m_0(x)\in L^1(\Om),\,\,m_0=0\,\,\mbox{if}\,\,\rho_0=0, \frac{|m_0|^2}{\rho_0}\in L^1(\Om),\\
H(x,0)=H_0(x)\in L^2(\Om),\quad \div H_0=0\,\,\mbox{in}\,\,\mathcal{D}'(\Om),\\
u|_{\partial \Om}=0,\quad H|_{\partial \Om}=0.
\end{cases}
\ee

 \begin{definition}\label{Def1}  
 A triple $(\rho,u,H)$ is called  a finite energy weak solution of the problem \eqref{MHD}, \eqref{IC}  if for any $T>0,$\\
$\bullet$~$\rho,u$ and $H$ belong to the following classes
$$ 0\le\rho\in L^{\infty}(0,\infty; L^\ga(\Om)),
$$ 
 $$ u\in L^2(0,\infty;(H_0^1(\Om))^3),\,\,
\rho|u|^2\in L^{\infty}(0,\infty; L^1(\Om)),
$$
$$
H \in L^\infty(0,\infty;L^2(\Om))\cap L^2(0,\infty; (H_0^1(\Om))^3),
 \,\div H=0.$$
$\bullet$  Eqs. \eqref{MHD} holds in $\mathcal{D}'(0,\infty; \Om).$ Moreover, Eqs. \eqref{MHD}$_1$  holds in $\mathcal{D}'(0,\infty;\R^3)$ provided $(\rho,u)$ was prolonged to be zero on $\R^3\setminus {\Om}.$  \\
$\bullet$~The energy $$
E(t):=\int_{\Om}\left(\frac12 \rho |u|^2+\frac{1}{\ga-1}\rho^{\ga}+\frac12|H|^2\right)dx
$$  is locally integrable on $(0,\infty)$ and  for  any $0\le \psi(t)\in \mathcal{D}(0,\infty),$ it holds
\be\la{qj2} -\int_0^\infty \psi_tE(t)dt
+\int_0^\infty\psi \int_\Om \left[\mu|\na u|^2+(\lambda+\mu)(\div u)^2+\nu|\na\times H|^2\right]dxdt\leq 0.
\ee \\ Moreover, \be\label{a.10}
\sup_{t\ge 0} E(t)+\int_0^t\int_{\Om}\left[\mu|\na u |^2+(\lambda+\mu)|\div u |^2+\nu |\nabla \times H |^2\right]dxdt\leq E_0
 \ee
  where
$$
E_0:=\int_\Om \left( \frac{|m_0|^2}{2\rho_0}
+\frac12|H_0|^2+\frac{1}{\ga-1}\rho_0^\ga\right)dx<\infty.
$$
\end{definition}

In \cite{HW}, Hu-Wang proved
 \begin{lemma}[\cite{HW}]\label{lma.1} Assume that $\Om\subset\R^3$ is a bounded domain with a boundary of class $C^{2+\kappa},\kappa>0,$ and $\ga>\frac{3}{2}$. Then for any given $T>0,$  the initial-boundary value problem \eqref{MHD} and \eqref{IC} has a finite energy weak solution $(\rho,u,H)$ on $\Om\times (0,T)$ satisfying \eqref{1} and \eqref{hq1}.
 Moreover, there exist
 a stationary state of velocity $u_s=0,$ and a stationary state of magnetic field
 $H_s=0$ such that, for $\n_s$ as in \eqref{1}, as $t\rightarrow \infty,$
 \bnn
\begin{cases} \rho(x,t)\rightarrow \rho_s,\,\,\mbox{strongly in}\,L^\ga(\Om);\\
 u(x,t)\rightarrow u_s=0,\,\,\mbox{strongly in}\,L^2(\Om);\\
 H(x,t)\rightarrow H_s=0,\,\,\mbox{strongly in}\,L^2(\Om).
\end{cases}
\enn

 \end{lemma}

Now we can modify slightly our method to prove the following:
\begin{theorem}\label{Th1} Assume that  the conditions of Lemma \ref{lma.1} hold.   Then  there exist positive constants $C_1$ and $C_2$ both depending only on $\Om,A, \ga,\mu,\lm,\nu,\te_0,\n_s,$ and $ E_0$ such that $(\rho,u,H)$,  the  finite energy weak solution  to  \eqref{MHD} and \eqref{IC}  whose existence is guaranteed by  Lemma \ref{lma.1}, satisfies the following decay property:
\be\la{av22}
\ds\int_{\Om}(\rho|u|^2+|H|^2+G(\rho,\n_s))dx\leq C_2\exp\{-C_1t\}\,\,\mbox{a.e.}\,\,t>0,
\ee
with
\bnn G(\rho,\n_s):=\rho\int_{\rho_s}^\rho \frac{h^\ga-\rho_s^\ga}{h^2}dh.
\enn
\end{theorem}

\pf It follows from \eqref{IC}$_4$ and Sobolev's inequality that for $\te\in (0,2\ga/3],$ \be \la{av21}\ba &\left|\int_\Om (\nabla \times H)\times H\cdot\mathcal{B}([\rho^\theta]_\eps-\overline{[\rho^\theta]_\eps})dx\right|\\& \le C\|H\|_{L^6(\Om)}\|\na H\|_{L^2(\Om)}\|\mathcal{B}([\rho^\theta]_\eps- \overline{[\rho^\theta]_\eps})\|_{L^{3}(\Om)}\\&\le C\|\na H\|_{L^2(\Om)}^2,  \ea \ee due to \be\la{av23}  \|H\|_{L^6(\Om)} \le C\|\na H\|_{L^2(\Om)}. \ee
With \eqref{av21} and \eqref{av23}   at hand, one can follow the proof of Theorem \ref{Th} and obtain \eqref{av22}. \thatsall


\begin{thebibliography}{999}
\bibitem{A} Adams, R., Fournier, J.  Sobolev spaces. Second edition. Academic Press, New York, 2003.

\bibitem{B} Bogovskii, M. E. Solution of some vector analysis problems connected with operators div and grad. Trudy Sem. S. L. Sobolev. {\bf 80} (1980), 5-40 (in Russian).
\bibitem{CH} Cabannes, H. Theoretical Magnetofluiddynamics. Academic Press, New York, 1970.
\bibitem{F} Feireisl, E. Dynamics of Viscous Compressible Fluids. Oxford University Press, Oxford, 2004.

\bibitem{FP1} Feireisl, H., Petezeltov\'{a}, H. Large-time behavior of solutions to the Navier-Stokes equations of compressible flow. Arch. Ration. Mech. Anal. {\bf 150} (1999), 77-96.

 \bibitem{FNP} Feireisl, E., Novotn\'{y}, A., Petezeltov\'{a}, H. On the existence of globally defined weak solutions to the Navier-Stokes equations.
        J. Math. Fluid Mech. {\bf 3} (2001), 358-392.
 \bibitem{FZZ} Fang, D. Y., Zi, R. Z., Zhang, T. Decay estimates for isentropic compressible Navier-Stokes equations in bounded domain. J. Math. Anal. Appl. {\bf 386} (2012), 939-947.
   \bibitem{GA} Galdi,  G. P.  An Introduction to the Mathematical Theory of the Navier-Stokes Equations. I. Springer-Verlag. New York, 1994.
 \bibitem{HLX} Huang, X. D., Li, J., Xin, Z. P. Global well-posedness of classical solutions with large oscillations and vacuum to the three-dimensional
     isentropic compressible Navier-Stokes equations. Comm. Pure Appl. Math.
     {\bf 65} (2012), 549-585.
 \bibitem{H} Hoff, D. Global  existence for 1D, compressible, isentropic Navier-Stokes equations with large initial data. Trans. Amer. Math. Soc. {\bf 303} (1987), 169-181.
 \bibitem{H3} Hoff, D. Global solutions of the Navier-Stokes equations for multidimensional compressible flow with discontinuous initial data. J. Differntial Equations. {\bf 120} (2013), 215-254.
 \bibitem{HZ}Hoff, D., Zumbrun, K. Multi-dimensional diffusion waves for the Navier-Stokes equations of compressible flow. Indiana Univ. Math. J. {\bf 44} (1995), 604-676.
     \bibitem{HW} Hu,  X.,  Wang, D. Global existence and large-time behavior of solutions to the threedimensional
equations of compressible magnetohydrodynamic flows. Arch. Ration. Mech.
Anal., 197 (2010), 203-238.

 \bibitem{JZ} Jiang, S., Zhang, P. On sperically symmetric solutions of the compressible isentropic Navier-Stokes equations. Comm. Math. Phys.
       {\bf 215} (2001), 559-581.
 \bibitem{KK} Kagei, T., Kobayashi, T. On large time behavior of solutions to the compressible Navier-Stokes equations in the half space in $\R^3$. Arch. Ration. Mech. Anal. {\bf 165} (2002), 89-159.

 \bibitem{KS1} Kazhikhov, A. V., Shelukhin, V. V. Unique global solution with respect to time of initial-boundary value problems for one-dimensional equations of a viscous gas. J. Appl. Math. Mech. {\bf41} (1977),  273-282; translated from Prikl. Mat. Meh. {\bf41} (1977), 282-291.
\bibitem{KL}Kulikovskiy, A. G., Lyubimov, G. A.  Magnetohydrodynamics. Addison-Wesley, Reading, 1965.

 \bibitem{KS2}Kobayashi, T., Shibata, Y. Decay estimates of solutions for the equations of motion of compressible viscous and heat-conductive gases in an exterior domain in $\R^3$. Comm. Math. Phys. {\bf 200} (1999), 621-659.
 \bibitem{LL} Laudau, L. D., Lifshitz, E. M.  Electrodynamics of Continuous Media, 2nd edn. Pergamon, New York, 1984.

  \bibitem{L1} Lions, P. L. Mathematical topics in fluid mechanics. Vol. 1. Incompressible models. Oxford University Press, New York, 1996.
   \bibitem{L2} Lions, P. L. Mathematical topics in fluid mechanics. Vol. 2. Compressible models. Oxford University Press, New York, 1998.

   \bibitem{MN1} Matsumura, A., Nishida, T. The initial value problem for the equations of motion of compressible viscous and heat-conductive fluids. Proc. Jpn. Acad. Ser. A. {\bf 55} (1979), 337-342.
    \bibitem{MN2} Matsumura, A., Nishida, T. The initial value problem for the equations of motion of viscous and heat-conductive gases. J. Math. Kyoto Univ. {\bf 20} (1980), 67-104.
     \bibitem{MN3} Matsumura, A., Nishida, T. Initial boundary value problems for the equations of motion of compressible viscous and heat-conductive fluids.  Comm. Math. Phys. {\bf 89} (1983), 445-464.
    \bibitem{Na} Nash, J. Le probl\`{e}me de Cauchy pour les  \'{e}quations diff\'{e}rentielles d’un fluide g\'{e}n\'{e}ral. Bull. Soc. Math. France. {\bf 90} (1962), 487-497.

   \bibitem{P} Padula, M. On the exponential stability of the rest state of a viscous compressible fluids. J. Math. Fluid Mech. {\bf 1} (1999), 62-77.

   \bibitem{PO} Ponce, G. Global existence of small solution to a class of nonlinear evolution equations. Nonlinear  Anal. {\bf9} (1985), 339-418.

    \bibitem{Ser} Serrin, J. On the uniqueness of compressible fluid motion. Arch. Ration. Mech. Anal. {\bf 3} (1959), 271-288.

\end{thebibliography}
\end{document}